\documentclass{amsart}

\usepackage{amssymb,amsfonts, latexsym,amsmath,hhline,array,longtable,hhline}
\usepackage{pdfsync,color, comment,colortbl, mathrsfs,stmaryrd,cite,graphicx}

\usepackage{ifpdf}
\ifpdf \usepackage[colorlinks=true, citecolor=blue, linkcolor=blue, urlcolor=blue]{hyperref} \fi

\newcommand{\cal}{\mathcal}

\newtheorem{formula}{}[section]
\newtheorem{definition}[formula]{Definition}
\newtheorem{corollary}[formula]{Corollary}
\newtheorem{remark}[formula]{Remark}
\newtheorem{lemma}[formula]{Lemma}
\newtheorem{theorem}[formula]{Theorem}

\def\thrm{\begin{theorem}}
\def\thrml#1{\begin{theorem}\label{#1}}
\def\ethrm{\end{theorem}}
\def\rmrk{\begin{remark}}
\def\rmrkl#1{\begin{remark}\label{#1}}
\def\ermrk{\end{remark}}
\def\dfntn{\begin{definition}}
\def\dfntnl#1{\begin{definition}\label{#1}}
\def\edfntn{\end{definition}}
\def\nmrt{\begin{enumerate}}
\def\enmrt{\end{enumerate}}
\def\tm#1{\item[{\rm (#1)}]}
\def\qtnl#1{\begin{equation}\label{#1}}
\def\eqtn{\end{equation}}
\def\lmm{\begin{lemma}}
\def\lmml#1{\begin{lemma}\label{#1}}
\def\elmm{\end{lemma}}
\def\crllr{\begin{corollary}}
\def\crllrl#1{\begin{corollary}\label{#1}}
\def\ecrllr{\end{corollary}}
\def\css{\begin{cases}}
\def\ecss{\end{cases}}
\def\prf{\begin{proof}}
\def\eprf{\end{proof}}

\def\cA{{\cal A}}

\def\cX{{\cal X}}

\def\cY{{\cal Y}}

\def\mC{{\mathbb C}}

\def\mR{{\mathbb R}}
\def\mS{{\mathbb S}}

\def\fX{{\mathfrak X}}

\DeclareMathOperator{\aut}{Aut}
\DeclareMathOperator{\alt}{Alt}

\DeclareMathOperator{\id}{id}
\DeclareMathOperator{\im}{im}

\DeclareMathOperator{\inv}{Inv}

\DeclareMathOperator{\iso}{Iso}
\DeclareMathOperator{\mat}{Mat}
\DeclareMathOperator{\orb}{Orb}

\DeclareMathOperator{\rk}{rk}

\DeclareMathOperator{\sym}{Sym}
\DeclareMathOperator{\WL}{WL}

\def\bone{{\bf 1}}
\def\grp#1{\langle {#1}\rangle}
\def\scp#1{\langle\hspace{1pt} {#1}\hspace{1pt}\rangle}

\def\phmb#1{{\phantom{x}\hspace{-2mm}^{#1}}}
\def\qaq{\quad\text{and}\quad}

\def\wt{\widetilde}

\begin{document}

%\begin{comment}
\title{On the  Weisfeiler-Leman dimension of some polyhedral graphs}
\author{Haiyan Li }
\address{School of Science, Hainan University, Haikou, China}
\email{lhy9694@163.com}
\author{Ilia Ponomarenko}
\address{School of Science, Hainan University, Haikou, China and Steklov Institute of Mathematics at St. Petersburg, Russia}
\email{inp@pdmi.ras.ru}
\author{Peter Zeman}
\address{Technical University of Denmark, Denmark}
\email{zeman.peter.sk@gmail.com}
\thanks{}
\date{}

\begin{abstract}
Let $m$ be a positive integer, $X$ a graph with vertex set~$\Omega$, and $\WL_m(X)$ the coloring of the  Cartesian $m$-power~$\Omega^m$, obtained by the $m$-dimensional Weisfeiler-Leman algorithm. The $\WL$-dimension of the graph~$X$ is defined to be the smallest $m$ for which the coloring $\WL_m(X)$ determines~$X$ up to isomorphism. It is known that the $\WL$-dimension of any planar graph is~$2$ or~$3$,  but no  planar graph of $\WL$-dimension~$3$ is known. We prove that the $\WL$-dimension of a  polyhedral (i.e., $3$-connected planar) graph~$X$ is at most~$2$ if the color classes of the coloring $\WL_2(X)$ are the orbits of the componentwise action of the group $\aut(X)$ on~$\Omega^2$.
\end{abstract}
%\end{comment}

\maketitle

\section{Introduction}

The \emph{WL-dimension} $\dim_{\scriptscriptstyle\WL}(X)$  of a graph~$X$ is a numerical invariant  defined  (implicitly) by M.~Grohe in \cite{Grohe2017} in connection with the Graph Isomorphism   Problem; the name refers to the Weisfeiler-Leman algorithm which constructs a canonical coloring of  the pairs of  vertices of~$X$. Roughly speaking, the WL-dimension is equal to the minimal positive integer $m$ such that the graph~$X$ can be identified up to isomorphism in a fragment of the first-order logic with counting quantifiers and~$m$ variables. In general,  $1\le \dim_{\scriptscriptstyle\WL}(X)\le n$, where~$n$ is the vertex number of~$X$. 

The graphs of WL-dimension~$1$ have completely been characterized in  \cite{Arvind2015} and independently in~\cite{Kiefer2015}. A complete description of the graphs of WL-dimension~$2$ seems to be hopeless (it would imply a description of all strongly regular graphs identified by the parameters up to isomorphism); the difficulties that arise on the way to obtaining such a description are clearly visible in a very interesting paper~\cite{Fuhlbr2018a}. It should be mentioned that there are graphs of arbitrary large WL-dimension, see~\cite{CaiFI1992}.

In the last decade, quite a lot of papers have appeared in which the WL-dimension is calculated (or estimated from above by some small constant) for graphs belonging to a certain class. Not being able to give a detailed review of these papers, we will focus here only on paper~\cite{Kiefer2017}, which is directly related to  what follows. It is proved there that  the WL-dimension of every planar graph is at most~$3$. It is quite interesting that neither~\cite{Kiefer2017} nor more recent paper~\cite{Kiefer2022} contain an example of a planar graph $X$ with $\dim_{\scriptscriptstyle\WL}(X)=3$. An initial motivation of the present paper was to understand whether such examples exist among the \emph{polyhedral} (i.e.,  planar and $3$-connected) graphs.

With each (polyhedral) graph $X$ with vertex set $\Omega$, one can associate a \emph{coherent  configuration} $\WL(X)$, which can be thought of as a partition~$S$  of the Cartesian square~$\Omega^2$ into the color classes of the canonical coloring obtained by the Weisfeiler-Leman algorithm applied to~$X$ (the exact definitions are given in Section~\ref{061022a}). One of the main results of~\cite{Kiefer2022} shows that if the edge set $E$ of the polyhedral graph $X$ is a class of the partition~$S$, then $X$ is edge transitive, i.e., $E\in\orb(G,\Omega^2)$, where $G=\aut(X)$. All edge transitive polyhedral graphs are known and one can verify that for each of them the coherent configuration $\WL(X)$ is schurian, which means that each class of $S$ belongs to $\orb(G,\Omega^2)$. 

A graph $X$ is said to be \emph{schurian} if the coherent configuration $\WL(X)$  is schurian. Note that almost all graphs are schurian, and all graphs with at most~$13$ vertices are schurian (see~\cite{Klin2016}). In particular, this covers all the Platonic graphs except for the dodecahedral graph, for which this can be verified by a straightforward computation. We believe that every polyhedral graph is schurian, and if this is true, then by Theorem~\ref{310722a} below, the WL-dimension of any polyhedral graph is at most~$2$. 

\thrml{310722a}
Let $X$ be a polyhedral graph. If  $X$ is schurian, then $\dim_{\scriptscriptstyle\WL}(X)\le 2$.
\ethrm

It was proved in~\cite{Fuhlbr2018a} that if $X$ is an arbitrary graph, then $\dim_{\scriptscriptstyle\WL}(X)\le 2$ if and only if the coherent configuration $\cX=\WL(X)$ is separable, i.e., every algebraic isomorphism from $\cX$ to another coherent configuration is induced by isomorphism. Assume that the graph $X$ is schurian and $G=\aut(X)$. Then $\cX$ coincides with the coherent configuration~$\inv(G)$ corresponding to the partition of $\Omega^2$ into the orbits of the group $G\le\sym(\Omega^2)$. If, in addition, $X$ is polyhedral, then $G$ is a  finite spherical group, i.e., a finite subgroup of the orthogonal group $O(3,\mR^3)$, see, e.g.,~\cite{Klavik2022}. It follows that, in this case, $\cX$ is the coherent configuration of a \emph{strongly spherical} group, by which we mean a spherical group which is the full automorphism group of a polyhedral graph. Thus Theorem~\ref{310722a} immediately follows from the theorem below.

\thrml{310722b}
Let $G$ be a strongly spherical permutation group. Then the coherent configuration $\inv(G)$ is separable.
\ethrm

The proof of Theorem~\ref{310722b} is given in Section~\ref{011122w}. The main idea of the proof is to use a  recent classification of finite spherical permutation groups~\cite{Klavik2022}. Based on  some sufficient conditions for a coherent configuration to be separable (Section~\ref{021122a}), we reduce the proof to the study of  coherent configurations of   spherical groups $G$ which, with the exception of a finitely many groups, are  isomorphic to $C_n$, $D_{2n}$, $C_n\times C_2$, and $D_{2n}\times C_2$. First, we prove in a more or less standard way that  the coherent configuration $\inv(G)$ is separable if~$G$ belongs to the infinite families (Theorem~\ref{221022x}). The same idea supplied by computer calculations works for some groups not belonging to these families (Theorem~\ref{221022y}).

In the remaining cases $G$ is isomorphic to $\sym(4)$, $\alt(5)$, or $\sym(4)\times C_2$. To manage these cases, we introduce in Section~\ref{290222a} a concept of spherical representation of a coherent configuration and prove that under some natural conditions, a coherent configuration admitting such a representation in the dimension~$3$ is separable (Theorem~\ref{070922v}). Again, using computer calculations we prove that in the considered cases, the coherent configurations $\inv(G)$ satisfy the above conditions (Theorem~\ref{221022z}) and hence are separable.

In order to make the paper as self-contained as possible, we present in Section~\ref{061022a} the necessary notation and facts from the theory of  coherent configurations. More details about this theory can be found in monograph~\cite{CP2019}.

\section{Coherent configurations}\label{061022a}

 \subsection{Notation}
Throughout the paper, $\Omega$ denotes a finite set. For $\Delta\subseteq \Omega$, the Cartesian product $\Delta\times\Delta$ and its diagonal are denoted by~$\bone_\Delta$ and $1_\Delta$, respectively. For  $s\subseteq\bone_\Omega$, we set $s^*=\{(\alpha,\beta): (\beta,\alpha)\in s\}$, and $\alpha s=\{\beta\in\Omega:\ (\alpha,\beta)\in s\}$ for all $\alpha\in\Omega$. The composition $r\cdot s$ of the relations $r,s\subseteq\bone_\Omega$ consists of all pairs $(\alpha,\beta)\in\bone_\Omega$ such that $(\alpha,\gamma)\in r$ and $ (\gamma,\beta)\in s$ for some $\gamma\in\Omega$. For any collection~$S$ of relations, we denote by $S^\cup$ the set of all unions of elements of~$S$, and consider~$S^\cup$ as a poset with respect to inclusion.

The set of classes of  an equivalence relation $e$ on~$\Omega$ is denoted by $\Omega/e$. For $\Delta\subseteq\Omega$, we set $\Delta/e=\Delta/e_\Delta$ where $e_\Delta=\bone_\Delta\cap e$. If the classes of $e_\Delta$ are singletons, $\Delta/e$ is identified with $\Delta$. 

We use standard notations $\sym(n)$, $\alt(n)$,  $C_n$, and $D_{2n}$ for the symmetric and alternating groups of degree~$n$, and  for cyclic and dihedral groups of orders~$n$ and~$2n$, respectively. 

The algebra of all square matrices with entries belonging to a field $F$ and rows and columns indexed by the elements of~$\Omega$  is denoted by $\mat_\Omega(F)$.  The points of~$\Omega$ are treated as vectors of the $F$-linear space spanned by $\Omega$. In particular, if $\alpha\in\Omega$ and $A\in\mat_\Omega(F)$, then $A\alpha$ is the $\alpha$-column of~$A$.

The Hadamard product $A\circ B$ of the matrices $A,B\in\mat_\Omega(\mR)$ is defined by the formula $(A\circ B)_{\alpha,\beta}=A_{\alpha,\beta}B_{\alpha,\beta}$ for all $\alpha,\beta\in\Omega$. The adjacency matrix of a relation $s\subseteq\bone_\Omega$ is denoted by $A_s$. Clearly, $A_{s^*}$ is equal to the transpose $(A_s)^*$ of the matrix~$A_s$. 

\subsection{Rainbows}
Let $S$ be a partition of $\Omega^2$. A pair $\mathcal{X}=(\Omega,S)$ is called a \emph{rainbow} on $\Omega$ if 
\nmrt
\tm{C1}  $1_\Omega\in S^\cup$,
\tm{C2} $s^*\in S$ for all $s\in S$.
\enmrt
The rainbow $\cX$ is said to be \emph{symmetric} if $s^*=s$ for all $s\in S$. The numbers $|\Omega|$ and $\rk(\cX)=|S|$ are called the {\it degree} and {\it rank} of~$\cX$. The elements of $S$ and of $S^\cup$ are called {\it basis relations} and \emph{relations} of~$\cX$. A unique basis relation containing the pair~$(\alpha,\beta)$ is denoted by~$r(\alpha,\beta)$.  The set of all relations is closed with respect to intersections and unions.

Let $\cX'=(\Omega', S')$ be a rainbow. A  {\it combinatorial isomorphism} or, briefly, \emph{isomorphism} from $\cX$ to $\cX'$   is defined to be a bijection $f: \Omega\rightarrow \Omega'$ such that  the relation $s^f=\{(\alpha^f, \beta^f): (\alpha, \beta)\in s\}$ belongs to~$S'$ for all $s\in S$. In this case,  the rainbows are said to be {\it isomorphic}; the set of all isomorphisms from $\cX$ to~$\cX'$ is denoted by $\iso(\cX,\cX')$. The group of all isomorphisms from $\cX$ to itself contains a normal subgroup
$$
\aut(\cX)=\{f\in\sym(\Omega):\ s^f=s \text{ for all } s\in S\}
$$
called the {\it automorphism group} of $\cX$.

There is a natural partial order\, $\le$\, on the set of all rainbows  on~$\Omega$. Namely, given two such rainbows~ $\cX$ and $\cY$, we set
$$
\cX\le\cY\ \Leftrightarrow\ S(\cX)^\cup\subseteq S(\cY)^\cup.
$$
The minimal and maximal elements with respect to this order are the {\it trivial} and {\it discrete} rainbows, respectively; in the former case, $S$ consists of~$1_\Omega$ and its complement (unless $\Omega$ is not a singleton), and in the last case, $S$ consists of singletons. 
\begin{comment}
Note that the functor $\cX\to\aut(\cX)$  reverse the inclusion, namely,
\qtnl{140322a}
\cX\le \cY \Rightarrow  \aut(\cX)\ge \aut(\cY).
\eqtn
\end{comment}

\subsection{Coherent configurations}
A rainbow $\cX$ is called a \emph{coherent configuration} if, in addition to conditions~(C1) and~(C2), it satisfies the condition
\nmrt
\tm{C3} given $r,s,t\in S$, the number $c_{r,s}^t=|\alpha r\cap \beta s^{*}|$ does not depend on $(\alpha,\beta)\in t$. 
\enmrt
We say that the $c_{rs}^t$ are \emph{intersection numbers} of~$\cX$ and extend the notation by admitting $r$ and $s$ to be arbitrary
relations of~$\cX$: $c_{r,s}^t=\sum_{r',s'}c_{r', s'}^t$, where $r',s'$ run over the basis relations contained in $r$ and $s$, respectively.
The coherent configuration is said to be a {\it scheme} if $1_\Omega\in S$,  and {\it commutative} if $c_{r,s}^t=c_{s,r}^t$ for all $r,s,t\in S$.  

Let $s\in S$. Then given a reflexive $t\in S$, the intersection number $n_s=c_{s,s^*}^t$  is either zero or is equal to $|\alpha s|$ for all $\alpha$ with $\alpha s\ne\varnothing$;  this number  is called the {\it valency} of~$s$.  We say that $s$ is {\it thin} if $n_{s^{}}=n_{s^*}=1$. A coherent configuration is said to be \emph{semiregular} if every its basis relation is thin, and \emph{quasi-thin} if $n_s\le 2$ for all $s\in S$; a semiregular scheme is said to be \emph{regular}. 

Any $\Delta\subseteq\Omega$ such that $1_\Delta\in S$ is called a {\it fiber} of $\cX$. In view of condition~(C1), the set  of all fibers forms a partition of~$\Omega$. Every basis relation is contained in the Cartesian product of two uniquely determined fibers. Any union  of fibers is called a {\it homogeneity set} of~$\cX$.  For any two homogeneous sets $\Delta$ and $\Gamma$, we denote by~$S_{\Delta,\Gamma}$ the set of all  basis relations contained in $\Delta\times\Gamma$, and abbreviate $S_\Delta=S_{\Delta,\Delta}$.

Let $G\le\sym(\Omega)$. Denote by $S$ the set of all orbits $(\alpha,\beta)^G$ in the componentwise action of~$G$ on~$\Omega\times\Omega$, where $\alpha,\beta\in \Omega$. Then the pair 
$$
\inv(G)=\inv(G,\Omega)=(\Omega,S)
$$ 
is a coherent configuration. Any coherent configuration $\cX$ associated with permutation group $G$ in this way is said to be {\it schurian}. In this case every  fiber $\Delta$ of~$\cX$  is an orbit of~$G$,  and the points of $\Delta$ can be identified with the right cosets of a point stabilizer $G_\delta$ with $\delta\in\Delta$,  so that the action of~$G$ on~$\Delta$ is induced by right multiplications. Note that a coherent configuration~$\cX$ is schurian if and only if $\cX=\inv(\aut(\cX)$.  A  coherent configuration $\inv(G,\Omega)$ is quasi-thin if $|G_\alpha|\le 2$ for all $\alpha\in\Omega$.

A bijection $\varphi:S\rightarrow S'$ is called an \emph{algebraic isomorphism} from $\cX$ to $\cX'$ if for all $r,s,t\in S$, we have
$$
c_{r^\varphi,s^\varphi}^{t^\varphi}=c_{r,s}^t.
$$
In this case, $|\Omega'|=|\Omega|$. The algebraic isomorphism $\varphi$ is extended in a natural way to a poset isomorphism $S^\cup\to (S')^\cup$, denoted also by~$\varphi$. Thus, $\varphi$ respects the set theoretical operations on~$S^\cup$, and  $1^{}_{\Omega'}=\varphi(1_{\Omega^{}})$.

Let $\Phi$ be a group of algebraic automorphisms of the coherent configuration~$\cX$. Denote by $S^\Phi$  the set of all relations $\cup_{\varphi\in\Phi}\varphi(s)$,  $s\in S$, Then the pair $\cX^\Phi=(\Omega,S^\Phi)$ is a coherent configuration called an \emph{algebraic fusion} of~$\cX$.

Every isomorphism $f:\cX\to\cX'$ \emph{induces} the algebraic isomorphism $\varphi_f:\cX\to \cX',$ $s\mapsto s^f$. Not every algebraic isomorphism is induced by isomorphism. The coherent configuration $\cX$  is said to be \emph{separable} if every algebraic isomorphism from $\cX$  is induced by  isomorphism.  All trivial and discrete coherent configurations are separable.

\subsection{Adjacency algebra} Let $\cX$ be a coherent configuration on~$\Omega$. Then the linear space $\cA(\cX)$ consisting of $\mC$-linear combinations of the matrices $A_s$, $s\in S$, is a subalgebra of the full matrix algebra $\mat_\Omega(\mC)$. The algebra $\cA(\cX)$, called the \emph{adjacency algebra} of~$\cX$,  contains the all ones matrix, and is closed with respect to transpose and the Hadamard multiplication. 

Every algebraic isomorphism $\varphi:\cX\to\cX'$  induces a matrix algebra isomorphism $\wt\varphi:\cA(\cX)\to\cA(\cX')$ such that $\varphi(A_s)=A_{\varphi(s)}$ for all $s\in S$. The isomorphism $\wt\varphi$ respects the transpose and the Hadamard product. In what follows, we simplify notation and write $\varphi$ instead of $\wt\varphi$. 

\lmml{051022a}
In the above notation, let $A\in\cA(\cX)$ and $A'=\varphi(A)$. Then
\nmrt
\tm{1} $\{A^{}_{\alpha^{},\beta^{}}:\ \alpha,\beta\in\Omega\}= \{A'_{\alpha',\beta'}:\ \alpha',\beta'\in\Omega'\}$,
\tm{2} $A$ is symmetric if and only if so is $A'$, 
\tm{3} $A$ and $A'$ have the same eigenvalues,
\tm{4} $A\alpha=A\beta$ for some $\alpha\ne\beta$ if and only if  $A'\alpha'=A'\beta'$ for some $\alpha'\ne\beta'$.
\enmrt
\elmm
\prf
(1) Denote by $W$ and $W'$ the sets  in the left- and right-hand sides of the required equality, respectively. For each $w\in W$, denote by $s_w$ the relation of $\cX$ consisting of all pairs $(\alpha,\beta)$ such that $(A)_{\alpha,\beta}=w$. Then $A_{s_w}\in\cA(\cX)$ and
\qtnl{030922g}
\sum_{w'\in W'}w'A_{s_{w'}}=A'=\varphi(A)=\sum_{w\in W}w\,\varphi(A_{s_w})=\sum_{w\in W}wA_{\varphi(s_w)},
\eqtn
which shows that $W=W'$, because  $A_{\varphi(s_w)}$ and $A_{s_{w'}}$are \{0,1\} matrices for all $w$ and $w'$.

(2) It is known that $\varphi(s^*)=\varphi(s)^*$ for all $s\in S$. It follows that $\varphi(A)^*=\varphi(A^*)$. Since $A$ is symmetric if and only if $A=A^*$, we are done.

(3) Since $\varphi$ is a matrix algebra isomorphism, the matrices $A$ and $A'$ have the same minimal polynomial. Thus $A$ and $A'$ have the same eigenvalues.

(4) The equality $A\alpha=A\beta$ means that $A_{\gamma,\alpha}=A_{\gamma,\beta}$ for every $\gamma\in\Omega$. The last equality holds if and only if   $(\gamma,\alpha),(\gamma,\beta)\in s_w$ for some $w\in W$. However, the number of all $\gamma\in\Omega$ for which $(\gamma,\alpha),(\gamma,\beta)\in s_w$ is equal to the intersection number $c_{(s_w)^*,s_w}^r$,
where $r=r(\alpha,\beta)$. Thus, 
$$
A\alpha=A\beta\ \Leftrightarrow\ \sum_{w\in W}c_{(s_w)^*,s_w}^r=|\Omega|.
$$
Since the right-hand side equality is preserved by $\varphi$, we are done.
\eprf

\subsection{Parabolics and quotients}\label{140322w}  A relation of $\cX$  that is an equivalence relation on $\Omega$ is called a \emph{parabolic} of~$\cX$; the set of all parabolics  is denoted by~$E=E(\cX)$.  Given $e\in E$, we put $S_{\Omega/e}=\{s_{\Omega/e}:\ s\in S\}$, and put $S_\Delta=\{s_\Delta:\ s\in S,\ s_\Delta\ne\varnothing\}$ for any $\Delta\in\Omega/e$. The pairs 
$$
\cX_{\Omega/e}=(\Omega/e,S_{\Omega/e})\qaq \cX_\Delta=(\Delta,S_\Delta)
$$  
are coherent configurations called a \emph{quotient} of~$\cX$ modulo~$e$ and  a \emph{restriction} of~$\cX$ to~$\Delta$, respectively.  Note that the notation $\cX_\Delta$ is correct for every homogeneity set~$\Delta$, because such $\Delta$ is a class of an obvious parabolic with two classes. 

Let $\varphi:\cX\to\cX'$ be an algebraic isomorphism and  $e\in E$. Then  $\varphi(e)$ is a parabolic of~$\cX'$ with the same number of classes. Moreover,  $\varphi$ induces the algebraic isomorphism
\qtnl{041122w}
\varphi^{}_{\Omega^{}/e^{}}:\cX^{}_{\Omega^{}/e^{}}\to \cX'_{\Omega'/e'},\ s^{}_{\Omega^{}/e^{}}\mapsto s'_{\Omega'/e'},
\eqtn
where $e'=\varphi(e)$ and $s'=\varphi(s)$.  

\subsection{Coherent closure}
The {\it coherent closure} $\WL(T)$ of a set $T$ of binary relations on $\Omega$, is defined to be the smallest coherent configuration on $\Omega$, for which each relation of~$T$ is the union of some basis relations. When $T$ consists of the arc set of a graph $X$ (respectively the basis relations of a rainbow~$\cX$), we denote $\WL(T)$ by $\WL(X)$ (respectively, $\WL(\cX)$). Every isomorphism from a rainbow $\cX$ to a rainbow~$\cX'$ is also an isomorphism from  $\WL(\cX)$ to $\WL(\cX')$.

\section{Separable coherent configurations}\label{021122a}
In this section, we recall and prove some facts about separable coherent configurations. The first of them  is just a reformulation of~\cite[Theorem 3.3.19]{CP2019}.

\begin{comment}
	First,  every regular scheme is separable. 
	
	Next, the direct sum (respectively, tensor product) of coherent configurtions is separable if and only if each of them is separable, see \cite[Corollary~3.2.8]{CP2019} (respectivly, \cite[Lemma~3.3.24]{CP2019}).  
	
	The regular schemes form a small subclass of a so called partly regular coherent configurations \cite[Subsec.~3.3.3]{CP2019}. All of them are schurian and the lemma below
\end{comment}

\lmml{181022a}
Let $G$ be a permutation group. Assume that at least one of the orbits of~$G$ is faithful and regular. Then the coherent configuration~$\inv(G)$ is separable.
\elmm

Let  $\Delta$ and $\Gamma$ be orbits of a group $G\le\sym(\Omega)$. We say that  $\Delta$ is \emph{dominated} by~$\Gamma$, or that $\Gamma$ \emph{dominates}~$\Delta$, if for any $\delta\in\Delta$ there is  $\gamma\in\Gamma$ such that $G_\gamma\le G_\delta$. In this case, $|\Delta|\le |\Gamma|$, and we write $\Delta\preceq\Gamma$ or $\Gamma\succeq\Delta$. Obviously, a singleton orbit is dominated by any other orbit, and a faithful regular orbit dominates any other orbit.

One can see that $\Gamma\succeq\Delta$ if and only if  there is $s\in\orb(G,\Gamma\times\Delta)$ such that $|\gamma s|=1$ for all $\gamma\in\Gamma$. It is also clear that if $\Delta$ and $\Gamma$ are \emph{equivariant} orbits of~$G$ (i.e., the permutation groups~$G^\Delta$ and $G^\Gamma$ are equivalent), then each of them is dominated by the other. Note that the relation ``to be dominated'' is transitive. The following lemma is a special case of  \cite[Lemma~3.3.20]{CP2019}.

\lmml{310722d}
Let $G\le\sym(\Omega)$ and $\Delta\in\orb(G)$. Assume that $\Delta$ is dominated by another orbit. Then  the coherent configuration $\inv(G)$ is separable if and only if so is the coherent configuration $\cX=\inv(G,\Omega\setminus\Delta)$. Moreover, $\aut(\cX)\cong \aut(\inv(G))$.
\elmm

The hypothesis of Lemma~\ref{310722d} is obviously satisfied if $\Delta$ is a singleton. Furthermore, assume that $\Delta\ne\Omega$ is of cardinality~$2$ and is not dominated by another orbit. Then $\alpha s=\Delta$ for all  $\alpha\not\in\Delta$ and all $s\in S$ such that $\alpha s\,\cap\,\Delta\ne\varnothing$. It follows that the coherent configuration $\inv(G)$ is the direct sum of the coherent configurations $\inv(G,\Delta)$ and $\inv(G,\Omega\setminus\Delta)$. Since also any homogeneous two-point coherent configuration is trivial and hence is separable, we conclude that $\inv(G,\Delta)$  is separable. Thus,  the coherent configuration $\inv(G)$ is separable  if and only if so is the coherent configuration $\inv(G,\Omega\setminus\Delta)$, see \cite[Corollary~3.2.8]{CP2019}. Thus we arrive to the following statement.

\crllrl{200922a}
Let $G\le\sym(\Omega)$ and $\Delta\in\orb(G)$. Assume that $|\Delta|=1$ or~$2$. Then the coherent configuration $\inv(G)$ is separable if and only if so is the coherent configuration $\inv(G,\Omega\setminus\Delta)$.
\begin{comment}
	The conclusion of Lemma~\ref{310722d} holds true in each of the following cases:
	\nmrt
	\tm{a} there is an orbit of $G$, which is different from $\Delta$ and equivariant to~$\Delta$,
	\tm{b} $G$ has a faithful regular orbit different from~$\Delta$,
	\tm{c} $|\Delta|=1$.
	\enmrt
\end{comment}
\ecrllr

The group $G$ is said to be \emph{domination free} if no orbit $\Delta$ of $G$ dominates an orbit other than $\Delta$.  In this case, no orbit of $G$ is faithful regular or a singleton (Corollary~\ref{200922a}), and no two distinct orbits of $G$ are equivariant.

\lmml{240922b}
Let $G\cong D_{2n}$ be a permutation group, $n\ge 2$. Assume that every orbit of~$G$ is faithful of cardinality~$n$.  Then $\inv(G)$ is separable.
\elmm
\prf
By Lemma~\ref{310722d}, we may assume that $G$ is domination free. Then $G$ has no equivariant orbits and hence~$G$ has one or two orbits depending on whether $n$ is odd or even. First, assume that~$G$ has exactly one orbit. Then $\inv(G)$ is a scheme, and even a quasi-thin one, because $n_s\le |G_\alpha|=2$ for all its points $\alpha$ and all its basis relations~$s$. Moreover, the number $n_1$ of thin basis relations of  $\inv(G)$ is equal to $1$ or $2$ depending on whether $n$  is even or odd. On the other hand, according to~\cite[Theorem~1.1]{MuzP2012b}, a quasi-thin scheme of degree~$n$ is separable if $n_1< 4$. Thus the scheme $\inv(G)$ is separable.

Now assume that $G$ has exactly two orbits, say $\Delta$ and $\Gamma$. Then they are faithful and non-equivariant, and  $n$ is even. Take points $\delta\in\Delta$ and $\gamma\in\Gamma$,  and identify the sets $\Delta$ and $\Gamma$ with the right cosets of the subgroups $G_\delta$ and $G_\gamma$ respectively. Each of these cosets contains a unique element belonging to the subgroup $C_n\le G$. The  mappings 
$$
\Delta\to\Gamma,\ G_\delta c\mapsto G_\gamma c \qaq \Gamma\to\Delta,\ G_\gamma c\mapsto G_\delta c ,
$$
where $c\in C_n$,  yield  a permutation $f\in N_{\sym(\Omega)}(G)$ interchanging  $\Delta$ and~$\Gamma$.  Denote by $S$ the set of basis relations of the coherent configuration $\cX=\inv(G)$.  A straightforward computation shows  that  $f$ induces an involution $\varphi\in\sym(S)$ that interchanges $S_\Delta$ and $S_\Gamma$, and such that 
$$
s^\varphi=s^*\quad\text{for all}\ \,s\in S_{\Delta,\Gamma}\cup S_{\Gamma,\Delta},
$$
in particular, $f$ and $\varphi$ are  respectively,  an isomorphism and an algebraic isomorphism from~$\cX$ to itself, and $\Phi=\grp{\varphi}$ is a group of order~$2$. 

The algebraic fusion $\cX^\Phi$ is a schurian scheme associated with permutation group generated by~$G$ and~$f$.  Every point stabilizer of the latter group is also a point stabilizer of~$G$, and hence  is of order~$2$.  Therefore the scheme $\cX^\Phi$ is quasi-thin. It  has exactly one thin basis relation, and using ~\cite[Theorem~1.1]{MuzP2012b} again, we conclude that $\cX^\Phi$ is separable. This implies   that so is $\cX$ by \cite[Theorem 3.1.29]{CP2019}.
\eprf

We complete the section by a very special sufficient condition for a coherent configuration to be separable. To formulate the corresponding statement, we need some preparation.

Let $\cX$ be a schurian coherent configuration on~$\Omega$  and $\sigma\in\aut(\cX)$ a central involution. Put $\Delta=\{\alpha\in\Omega:\ \alpha^\sigma\ne\alpha\}$ and denote by $e=e(\sigma)$ the equivalence relation with classes $\{\alpha,\alpha^\sigma\}$, $\alpha\in\Omega$. Then $\Delta$ is a homogeneity set of~$\cX$ and $e$ is a parabolic of~$\cX$ such that $e_{\Omega\setminus\Delta}=1_{\Omega\setminus\Delta}$. The surjection 
\qtnl{231022w}
\pi:\Omega\to\Omega/e,\quad \alpha\mapsto\alpha e,
\eqtn
defines the quotient coherent configuration $\cX_{\Omega/e}=(\Omega/e, S_{\Omega/e})$ of $\cX$ modulo~$e$, where $S_{\Omega/e}=\{\pi(s):\ s\in S\}$. 

We say that a  parabolic $e^\bot$ of the coherent configuration $\cX_\Delta$ is  a \emph{complement} of~$e$ if any class of $e^\bot$ intersects any class of $e_\Delta$ in exactly one point; in this case the  parabolic $e^\bot$ has exactly two classes. If the complement $e^\bot$ exists, then for every $s\in S$, we have
\qtnl{300922w1}
s=\css
1_{\Delta(s)}\text{ or }   e_{\Delta(s)}\setminus   1_{\Delta(s)}  &\text{if $s\subseteq e_\Delta$,}\\
s_0\setminus e^\bot \text{ or }   s_0\cap e^\bot     &\text{if $s\in \bone_\Delta\setminus 1_\Delta$,}\\
s_0  &\text{otherwise.}\\
\ecss
\eqtn
where $\Delta(s)$ is the support of~$s$ and $s_0=\pi^{-1}(\pi(s))$. Moreover, one can see that if $\Delta=\Omega$, then  $\cX$ is isomorphic  to the tensor product of $\cX_{\Delta/e^{}}$ and $\cX_{\Delta/e^\bot}$, see \cite[Theorem~6]{Chen2021}, and the lemma below is a special case of \cite[Corollary~3.2.24]{CP2019},

\lmml{231022a}
Let $\cX$ be a schurian coherent configuration and $\sigma\in\aut(\cX)$ a central involution. Assume that the parabolic $e=e(\sigma)$ has a complement and the coherent configuration $\cX_{\Omega/e}$ is separable. Then $\cX$ is also separable.
\elmm
\prf
Let $\varphi:\cX\to\cX'$ be an algebraic isomorphism. Then $e'=\varphi(e)$ is a parabolic of $\cX'$ and $\varphi$ induces the algebraic isomorphism $\bar \varphi=\varphi_{\Omega/e}$, see~\eqref{041122w}. Since the coherent configuration $\cX_{\Omega/e}$ is separable, there is a bijection 
$$
\bar f:\Omega/e\to\Omega'/e'
$$
inducing $\bar\varphi$, where $\Omega'$ is the point set of~$\cX'$. Put $\Delta'=\Delta^\varphi$ and denote by $e^\bot$  a complement to~$e$. Then  $e'{\phmb{\bot}}=\varphi(e^\bot)$ is a  complement of the parabolic~$e'$ in the coherent configuration $\cX'_{\Delta'}$. Let $\Delta_1^{}$ and $\Delta_2^{}$ be the classes of~$e^\bot$, and let $\Delta'_1$ and~$\Delta'_2$ be the classes of~$e'{\phmb{\bot}}$.  Then $\Delta=\Delta^{}_1\cup\Delta^{}_2$ and $\Delta'=\Delta'_1\cup\Delta'_2$.

For any point $\alpha\in \Omega$, we define a point~$\alpha'\in\Omega'$ as follows. If $\alpha\not\in\Delta$, then $\{\alpha\}$ is a class of~$e$ and $\{\alpha\}^{\bar f}=\{\alpha'\}$ for some $\alpha'\not\in\Delta$. Now, let $\alpha\in\Delta$. Then there exist $i\in\{1,2\}$  and $\Lambda\in\Omega/e$ such that $\{\alpha\}= \Delta_i\cap \Lambda$.   Denote by $\alpha'$ the point of $\Delta'$ for which $\{\alpha'\}= \Delta'_i\cap \Lambda'$ with $\Lambda'=\Lambda^{\bar f}$. Now the mapping
$$
f:\Omega\to\Omega',\ \alpha\mapsto \alpha' 
$$
is a bijection taking  $\Delta^{}_1$ to $\Delta'_1$,  $\Delta^{}_2$ to $\Delta'_2$, and the complement of~$\Delta$ to the complement of~$\Delta'$. Let us verify that $f$ induces~$\varphi$. 

From the definition of~$f$, it follows that  $e^f=e'=\varphi(e)$ and $(e^\bot)^f=e'{\phmb{\bot}}=\varphi(e^\bot)$. In view of formula~\eqref{300922w1}, it therefore suffices to prove that for any $s\in S$,
$$
(1_{\Delta(s)})^f=\varphi(1_{\Delta(s)})%,\quad (\bone_{\Delta(s)})^f=\varphi(\bone_{\Delta(s)}),\quad
\qaq (s_0)^f=\varphi(s)_0,
$$
where $\varphi(s)_0=\pi'^{-1}(\pi'(\varphi(s))$. However,  $\varphi(1_{\Delta(s)})=1_{\Delta(\varphi(s))}$ and $\Delta(\varphi(s))=\Delta(s^f)$. Thus,
$$
(1_{\Delta(s)})^f=1_{\Delta(s^f)}=1_{\Delta(\varphi(s))}=\varphi(1_{\Delta(s)}),
$$
which proves the first equality. %The second one is proved similarly. 
Finally,  $\pi^{-1}(\pi(s))^f={\pi'}^{-1}(\pi'(s^f))$, because $e^f=e'$.  Since $f^{\Omega/e}=\bar f$ and $\bar f$ induces $\bar\varphi$, we obtain
$$
(s_0)^f=\pi^{-1}(\pi(s))^f={\pi'}^{-1}(\pi'(s^f))={\pi'}^{-1}(\pi(s)^{\bar f}))=
$$
$$
{\pi'}^{-1}(\bar\varphi(\pi(s^f)))={\pi'}^{-1}(\pi'(\varphi(s)))=\varphi(s)_0,
$$
as required.
\eprf

\section{Spherical representations of coherent configurations}\label{290222a}

Let $d\ge 1$ be an integer and  $\Omega$ a finite set.  Let  $\rho:\Omega\to\mR^{d+1}$ be an injective mapping  taking a point $\alpha\in\Omega$ to a point~$\rho(\alpha)$ of the real $d$-dimensional  unit sphere~$\mS_d\subset \mR^{d+1}$ centered at~$0$. Denote by $A_\rho\in\mat_\Omega(\mR)$ the Gram matrix of the vectors $\rho(\alpha)$, 
$$
(A_\rho)_{\alpha,\beta}=\scp{\rho(\alpha),\,\rho(\beta)},\quad \alpha,\beta\in\Omega,
$$
where  $\scp{.\,,\,.}$ is the scalar product in $\mR^{d+1}$. Thus $A_\rho$ is a symmetric matrix with unit diagonal.  

We say that the mapping  $\rho$ is  a $d$-dimensional \emph{spherical representation}, or, briefly, an {\it $\mS_d$-representation}  of  a coherent configuration~$\cX$ on~$\Omega$, if the matrix $A_\rho$ belongs to the adjacency algebra $\cA(\cX)$. The following lemma gives a simple sufficient condition for $\cX$ to have an $\mS_d$-representation. In what follows, $\| \cdot\|$ denotes the standard Euclidean norm.

\lmml{030922a}
Let $P\in\cA(\cX)$ be a real symmetric idempotent matrix of rank $d+1$ without zero columns and such that for all  $\alpha,\beta\in\Omega$,
\qtnl{010922a}
\alpha\ne\beta\quad\Rightarrow\quad
{{P\alpha}\over{\| P\alpha\|}}\ne{{P\beta}\over{\| P\beta\|}}.
\eqtn
Then the mapping $\rho_P:\alpha\mapsto {{P\alpha}\over{\| P\alpha\|}}$, $\alpha\in\Omega$, is an $\mS_d$-representation of~$\cX$. 
\elmm
\prf
It suffices to verify that $A_\rho\in\cA(\cX)$, where $\rho=\rho_P$.  Since $P$ is symmetric and $P^2=P$, we have 
$$
P_{\alpha,\beta}=\| P\alpha\|\,\| P\beta\|\,\scp{\rho(\alpha),\,\rho(\beta)}
$$ 
for all $\alpha,\beta\in\Omega$. It follows that $A_\rho=DPD$, where $D$ is the diagonal matrix  with entries $D_{\alpha,\alpha}=\frac{1}{\| P\alpha\|}$, $\alpha\in\Omega$. Note that  $D=I\circ P$, where $I$ is the identity matrix. Since $I$ and $P$ belong to~$\cA(\cX)$, this algebra contains $D$, and hence  contains $A_\rho=DPD$, as required.
\eprf

Every coherent configuration has a trivial $(n-1)$-dimensional spherical representation with $n=|\Omega|$; in this representation, $\rho(\alpha)$ is a \{0,1\} vector  the only nonzero entry of which is on the place~$\alpha$. Of particular interest are low dimensional spherical representations. 

\crllrl{070922a}
The coherent configuration of a polyhedral graph admits an $\mS_d$-representation with $d\le 2$.
\ecrllr
\prf
Let $\cX$ be the coherent configuration of  a polyhedral graph with Laplacian~$L$. Denote by $P$ the orthogonal projection to the eigenspace of the matrix~$L$, associated with the second minimal eigenvalue.  Then $P\in\mat_\Omega(\mR)$ is a  symmetric idempotent matrix. Moreover, the rank of~$P$ is at most~$3$ and $P$ has no zero columns, see Corollary~13.10.2 and Lemma~13.10.1 in~\cite{Godsil2001}. In particular,  condition~\eqref{010922a} is satisfied. 

By  Lemma~\ref{030922a}, it remains to verify that $P\in\cA(\cX)$. However, $P=p(L)$, where $p(\cdot)$ is a polynomial with real coefficients, see e.g.~\cite[Sec.~8.12]{Godsil2001}.  Since $L\in\cA(\cX)$, this implies that $P\in\cA(\cX)$, as required.\eprf

Let $\rho$ be an  $\mS_d$-representation of a coherent configuration~$\cX$, and let $W:=W(\rho)$ be the set of all entries of the matrix~$A_\rho$. Clearly, $1\in  W$. Denote by $S_\rho$ the set of all binary relations
\qtnl{301022s}
s(w):=s_\rho(w)=\{(\alpha,\beta)\in\Omega^2:\ (A_\rho)_{\alpha,\beta}=w \},
\eqtn
where $w\in W$. 

Obviously, $s(1)=1_\Omega$ and $s(w)^*=s(w)$ for all~$w$. Consequently,  $\cX_\rho=(\Omega,S_\rho)$ is a rainbow.   Moreover, the adjacency matrix of the relation $s(w)$  belongs to the algebra $\cA(\cX)$, because $A_\rho\in\cA(\cX)$ and $\cA(\cX)$ is closed with respect to the Hadamard multiplication. It follows that
\qtnl{170223f}
\cX\ge \cX_\rho,
\eqtn
in particular, every basis relation of the coherent configuration~$\cX$ is contained in a unique basis  relation of the rainbow~$\cX_\rho$.  Given $s\in S$, we denote by $w(s)$  a uniquely determined element $w\in W$ such that $s\subseteq s(w)$.

The inclusion \eqref{170223f} implies that both $s(1)$ and $s(-1)$ are relations of~$\cX$ (the second of them can be empty). Furthermore, the union $e=s(1)\cup s(-1)$ is an equivalence relation on $\Omega$ with classes of cardinality at most two. In particular, $e$ is a parabolic of $\cX$;  we call it the \emph{antipodal parabolic} with respect to $\rho$. A pair $(\alpha,\beta)$ belongs to~$e$ if and only if $\rho(\alpha)=\pm\rho(\beta)$; any other pair is said to be \emph{non-antipodal}.

\lmml{030922h}
Let $\varphi:\cX\to \cX'$ be an algebraic isomorphism. Assume that $\cX$ has an $\mS_d$-representation~$\rho$. Then $\cX'$ has an $\mS_d$-representation~$\rho'$ such that $W(\rho')=W(\rho)$ and $w(\varphi(s))=w(s)$ for all $s\in S$. %In particular, $\varphi$ takes $\cX^{}_{\rho^{}}$ to~$\cX'_{\rho'}$.
\elmm
\prf
By statements (2) and (3) of Lemma~\ref{051022a}, the matrices $A_\rho$ and $A'=\varphi(A_\rho)$ are symmetric and have the same real eigenvalues. These are nonnegative, because $A_\rho$ is the Gram matrix. It follows that $A'$ is the Gram matrix of some vectors $x(\alpha')\in\mR^{\Omega'}$, $\alpha'\in\Omega'$, such that
$$
(A')_{\alpha',\beta'}=\scp{x(\alpha'),\,x(\beta')},\quad \alpha',\beta'\in\Omega'.
$$

By statement (1)  of Lemma~\ref{051022a},  the set $W'$ of all entries of the matrix~$A'$ is equal to $W(\rho)$. Furthermore, the diagonal entries of $A'$ (coinciding with diagonal entries of $A$) are ones. Together with statement (4)  of Lemma~\ref{051022a},  this shows that the vectors $x(\alpha')$ are pairwise distinct points of a sphere $\mS_{d'}$ for some $d'$. Moreover, 
$$
d'+1=\rk(A')=\rk(A)=d+1.
$$ 
Thus, the mapping $\rho':\alpha'\mapsto x(\alpha')$ is an  $\mS_d$-representation of the coherent configuration of~$\cX'$ and $W(\rho)=W'=W(\rho')$.  Finally,  the equality $\varphi(s_\rho(w))=s_{\rho'}(w)$ easily follows from formula~\eqref{030922g}.
\eprf

In the sequel, we are interested in a sufficient condition for an $\mS_2$-representation~$\rho$ to be uniquely determined by some two distinct vectors in $\im(\rho)$. To this end, let $u,v\in W$ and $x,y\in \mS_2$. Put
$$
\Omega_{u v}(x, y)=\{z\in \mS_2:\ \scp{x,\,z}=u,\ \scp{y,\,z}=v\}.
$$

Assume that $\Omega_{u v}(x, y)\ne\varnothing$. The equations $\scp{x,\,z}=u$ and $ \scp{y,\,z}=v$ (both with respect to~$z$) define the circles on $\mS_2$: one is of radius $\| u-x\|$  centered in $x$ and the other is of radius $\| v-y\|$  centered in $y$. These two circles can coincide only if $x=\pm y$. Assume that these circles are distinct. Then the coordinates of the vector $z$ can easily be found by solving an appropriate square equation. Thus we obtain the following statement.

\lmml{040922a}
In the above notation, assume that $x\ne\pm y$. Then $|\Omega_{u v}(x, y)|\le 2$. Moreover, if $\pi$ is an isometry of the sphere~$\mS_2$, then
$$
\Omega_{u v}(x, y)^\pi=\Omega_{u v}(x^\pi, y^\pi).
$$
\elmm

In some cases, the set $\Omega_{u v}(x, y)$ can be interpreted in terms of the original coherent configuration $\cX$ equipped with an $\mS_2$-representation~$\rho$.  Namely, for any $u,v\in W$ and any $\alpha,\beta\in\Omega$, we put
$$
\Omega_{u v}(\alpha,\beta)=\alpha s(u)\,\cap\,\beta s(v).
$$
Note that if a point $\gamma$ belongs to this set, then  $ \scp{\rho(\alpha),\,\rho(\gamma)}=w(r(\alpha,\gamma))=u$ and $ \scp{\rho(\beta),\,\rho(\gamma)}=w(r(\beta,\gamma))=v$.  Therefore, $\rho(\gamma)\in \Omega_{u v}(\rho(\alpha),\rho(\beta)$. This shows that
\qtnl{040922c}
\rho(\Omega_{u v}(\alpha,\beta))\subseteq \Omega_{u v}(\rho(\alpha),\rho(\beta)).
\eqtn
Together with Lemma~\ref{040922a} this proves the following statement.

\crllrl{050922a1}
Let $u,v\in W$, and $t\in S$. Assume that $t$ is not contained in the antipodal parabolic of $\cX$ with respect to~$\rho$. Then $c_{s(u),s(v)}^t\le 2$ and  the equality is attained only if 
\qtnl{170423a}
\rho(\Omega_{u v}(\alpha,\beta))=\Omega_{u v}(\rho(\alpha),\rho(\beta))
\eqtn
tor all $(\alpha,\beta)\in t$.
\ecrllr

For any basis relation $s$ of $\cX$, contained in the antipodal parabolic, we obviously have $n_s=n_{s^*}=1$, i.e., $s$ is a thin relation. By Corollary~\ref{050922a1}, this implies the following statement.

\crllrl{070922a1}
Let a coherent configuration $\cX$ have an $\mS_2$-representation, and let $r,s,t\in S$. Assume that $t$ is not thin.  Then  $c_{r,s}^t\le 2$.
\ecrllr

Let a set $\Delta\subseteq\Omega$ consist of at least three points.  Given distinct $\alpha,\beta\in \Delta$, we denote by $S(\Delta;\alpha,\beta)=S_\rho(\Delta;\alpha,\beta)$ the set of all relations $t\in S$ contained in the intersection
$$
\bigcap_{\delta\in\Delta\setminus\{\alpha,\beta\}}r_\rho(\alpha,\delta)\cdot r_\rho(\delta,\beta)
$$
where $r_\rho(\alpha,\delta)$ and $r_\rho(\delta,\beta)$ are the basis relations of the rainbow $\cX_\rho$, containing the pairs $(\alpha,\delta)$ and $(\delta,\beta)$, respectively. We also define  $S(\Delta;\alpha,\beta):=\{r(\alpha,\beta)\}$ for two-element subsets $\Delta$ of~$\Omega$. Thus, in any case, $S(\Delta;\alpha,\beta)$ contains $r(\alpha,\beta)$. A meaning of the set $S(\Delta;\alpha,\beta)$ lies in the fact that, roughly speaking, no two relations in it can be distinguished using only the points from $\Delta$ and basis relations of the rainbow~$\cX_\rho$.
 
A key point in what follows is a concept of a $\rho$-closed set. Informally, together with any two non-antipodal points $\alpha$ and $\beta$ satisfying condition~\eqref{170423a}, such a set should contain each $\gamma\in \Omega_{u v}(\alpha,\beta)$.  Using  Corollary~\ref{050922a1}, we will show that this is true if 
 \qtnl{031022a}
 t\in S(\Delta;\alpha,\beta)\ \Rightarrow\ c_{s(u),s(v)}^t=2,
 \eqtn
where  $u=w(r(\alpha,\gamma))$ and $v=w(r(\beta,\gamma))$. 

To be more precise, the set $\Delta$ is said to be \emph{$\rho$-closed} if it contains every point $\gamma\in\Omega$ for which there exists a non-antipodal pair $(\alpha,\beta)\in\Delta^2$ such that the condition~\eqref{031022a} is satisfied.  In this case the right-hand side is true for $t=r(\alpha,\beta)$. A subset of~$\Omega$ is  said to be \emph{$\rho$-rigid} if the only $\rho$-closed set containing  it coincides with $\Omega$.  Finally, the representation~$\rho$ is said to be \emph{rigid} if there exists a $\rho$-rigid set of cardinality~$2$.

\lmml{070922b}
Let  a coherent configuration $\cX$ have an $\mS_2$-representation $\rho$. Assume that $(\mu,\nu)$ is a non-antipodal pair  forming  a $\rho$-rigid set. Then for any algebraic isomorphism $\varphi:\cX\to\cX'$ and any points $\mu',\nu'$ for which $r(\mu',\nu')=\varphi(r(\mu,\nu))$, there exists a bijection $f:\Omega\to\Omega'$ such that
\qtnl{031022a1}
(\mu,\nu)^f=(\mu',\nu') \qaq s_\rho(w)^f=\varphi(s_\rho(w))\quad\text{for all}\  w\in W(\rho).
\eqtn
\elmm
\prf
By Lemma~\ref{030922h},  the coherent configuration $\cX'$ has an $\mS_2$-representation~$\rho'$  such that $W(\rho)=W(\rho')$ and $w(s)=w(\varphi(s))$ for all $s\in S$. It follows that $w(r(\mu,\nu))=w(\varphi(r(\mu,\nu)))=w(r(\mu',\nu'))$ and hence $\scp{\rho(\mu),\,\rho(\nu)}=\scp{\rho'(\mu'),\,\rho'(\nu')}$. Therefore there is an isometry $\pi:\mS_2\to\mS_2$  such that 
$$
\rho(\mu)^\pi=\rho'(\mu')\qaq \rho(\nu)^\pi=\rho'(\nu').
$$
Denote by $\Delta$ the set of all $\gamma\in\Omega$ such that $\rho(\gamma)^\pi=\rho'(\gamma')$ for some $\gamma'\in\Omega'$. Clearly, $\mu,\nu\in\Delta$. \medskip

{\bf Claim.} {\it $\Delta=\Omega$.}\medskip

\prf
Assume on the contrary that $\Delta\ne\Omega$. Since $\{\mu,\nu\}$ is s $\rho$-rigid subset of~$\Delta$, the set $\Delta$ is not $\rho$-closed by the choice of $\mu$ and $\nu$. It follows that there are $\gamma\not\in\Delta$ and a non-antipodal pair $(\alpha,\beta)\in\Delta^2$ such that the implication~\eqref{031022a} is true. In particular $c_{s(u),s(v)}^t=2$, where  $t=r(\alpha,\beta) $. By Corollary~\ref{050922a1}, we obtain
$$
\rho(\Omega_{u v}(\alpha,\beta))=\Omega_{u v}(\rho(\alpha),\rho(\beta)).
$$
Since $u=w(r(\alpha,\gamma))$ and $v=w(r(\beta,\gamma))$, we have  $\gamma\in \Omega_{u v}(\alpha,\beta)$. Therefore,
$\rho(\gamma)\in \Omega_{u v}(\rho(\alpha),\rho(\beta))$. By Lemma~\ref{040922a}, 
$$
\rho(\gamma)^\pi\in \Omega_{u v}(\rho(\alpha),\rho(\beta))^\pi=\Omega_{u v}(\rho(\alpha)^\pi,\rho(\beta)^\pi).
$$
Since $\alpha,\beta\in\Delta$, there are $\alpha',\beta'\in\Omega'$ such that
$\rho'(\alpha')=\rho(\alpha)^\pi$ and $\rho'(\beta')=\rho(\beta)^\pi$.  This shows that 
\qtnl{100922t}
\rho(\gamma)^\pi\in \Omega_{u v}(\rho'(\alpha'),\rho'(\beta')).
\eqtn

The isometry $\pi$  yields the bijection from~$\Delta$ to  the set $\Delta'={\rho'}^{-1}(\rho(\Delta)^\pi)$, that  takes a point $\delta$ to the point ${\rho'}^{-1}(\rho(\delta)^\pi)$. This bijection is obviously an isomorphism from the rainbow $(\cX^{}_{\rho^{}})_{\Delta^{}}$ to the rainbow $(\cX'_{\rho'})_{\Delta'}$.  Since the relation  $t'=r(\alpha',\beta')$  belongs to $S(\Delta';\alpha',\beta')$, we conclude that 
$
\varphi^{-1}(t') \in S(\Delta;\alpha,\beta).
$
Now the  implication~\eqref{031022a} yields  $c_{s(u),s(v)}^{t''}=2$, where $ t''=\varphi^{-1}(t') $. Then by Corollary~\ref{050922a1}, we have
$$
\rho'(\Omega_{u v}(\alpha',\beta'))=\Omega_{u v}(\rho'(\alpha'),\rho'(\beta')).
$$
But then $\rho(\gamma)\in\rho'(\Omega(u,v;\alpha',\beta'))\subseteq \Omega'$, which contradicts the maximality of $\Delta$.
\eprf

By the claim, there exists a bijection $f:\Omega\to\Omega'$ such that $\rho'(\alpha^f)=\rho(\alpha)^\pi$. Obviously $(\mu,\nu)^f=(\mu',\nu')$ and
$$
w(r(\alpha,\beta)^f)=\scp{\rho'(\alpha^f),\,\rho'(\beta^f)}=\scp{\rho(\alpha)^\pi,\,\rho(\beta)^\pi}=
\scp{\rho(\alpha),\,\rho(\beta)}=w(r(\alpha,\beta)),
$$
as required.
\eprf

To state the main result of this section, we need one more definition. An $\mS_2$-representation of a coherent configuration~$\cX$ is said to be \emph{faithful} if $\WL(\cX_\rho)=\cX$. In all calculated examples, the coherent configuration of a polyhedral graph $X$ has a faithful $\mS_2$-representation. However, not every $\mS_2$-representation of an arbitrary graph~$X$ is faithful, e.g., if $\cX$ is the coherent configuration of the  M\"obius--Kantor graph, then the $\mS_2$-representations associated (as in the proof of Corollary~\ref{070922a}) with eigenvalues $-1$ and $1$ are not faithful.

\thrml{070922v}
Let  a coherent configuration $\cX$ have an $\mS_2$-representation~$\rho$ . Assume that $\rho$ is rigid and faithful. Then $\cX$ is schurian and separable.
\ethrm
\prf
Let $\varphi:\cX\to\cX'$ be an algebraic isomorphism. By Lemma~\ref{070922b}, there exists a rainbow isomorphism $f:\cX^{}_{\rho^{}}\to\cX'_{\rho'}$ such that the second part of formula~\eqref{031022a1} holds. Denote by $\psi$  the algebraic isomorphism induced by~$f$,
$$
\psi:\WL(\cX^{}_{\rho^{}})\to\WL(\cX'_{\rho'}),\ s\mapsto s^f.
$$ 
Since $\rho$ is faithful, we have $\WL(\cX^{}_{\rho^{}})=\cX$ and hence $\rk(\WL(\cX'_{\rho'}))=\rk(\cX')$, which implies $\WL(\cX'_{\rho'})=\cX'$, because $\cX'\ge \WL(\cX'_{\rho'})$.  Thus, $\psi$ is an algebraic isomorphism from $\cX$ to $\cX'$. 

Now if $\psi=\varphi$, then $\varphi$ is induced by the bijection $f$. Since $\varphi$ is arbitrary, the coherent configuration $\cX$ is separable. Assume that $\psi\ne\varphi$. Then $\Phi=\grp{\varphi\psi^{-1}}$ is a nontrivial group of algebraic automorphisms of~$\cX$. Hence, the algebraic fusion $\cX^\Phi$ is strictly smaller than $\cX$. Since also $\Phi$ leaves each basis relation of the rainbow $\cX_\rho$ fixed, we obtain
$$
\cX>\cX^\Phi=\WL(\cX^\Phi)\ge \WL(\cX_\rho)=\cX,
$$
a contradiction.

To prove the schurity of~$\cX$, it suffices to verify that given $(\mu,\nu),(\mu',\nu')\in s\in S$, one can find $f\in\aut(\cX)$ such that $(\mu,\nu)^f=(\mu',\nu')$. Assume first that $s$ is not contained in the antipodal parabolic~$e$ (with respect to~$\rho$). Let $\varphi=\id_\cX$ be the identical algebraic automorphism of~$\cX$. Then by Lemma~\ref{070922b}, there exists a bijection $f:\Omega\to\Omega$  such that both parts of formula~\eqref{031022a1} hold.  By the first part of the proof for $\cX=\cX'$, $f$ induces $\varphi=\id_\cX$ and hence $f\in\aut(\cX)$, as required. 

To complete the proof, assume that $s\subseteq e$. Without loss of generality, we may assume that $\Omega$ contains a point~$\lambda$  different from~$\nu$. Then  $r(\mu,\lambda)\not\subseteq e$ (recall that each class of $e$ has at most two points).  Since $s$ is the basis relation containing both $(\mu,\nu)$ and~$(\mu',\nu')$, there exists a point $\lambda'$ such that $r(\mu,\lambda)=r(\mu',\lambda')$. By the above paragraph, one can find $f\in\aut(\cX)$ such that $(\mu,\lambda)^f=(\mu',\lambda')$.  Now if $\lambda=\mu$, then $\lambda'=\mu'$ and we are done. Assuming $\lambda\ne\mu$, we have $\lambda'\ne\mu'$. On the other hand, $e^f=e$, because $e$ is a parabolic of~$\cX$. Consequently, $\mu'=\mu^f$ and $\nu^f$ are in the same class of~$e$. Since the class containing $\mu'$ contains $\nu'$, we conclude that $\nu^f=\nu'$, as required.  
\eprf

\section{Proof of Theorem~\ref{310722b}}\label{011122w}

\subsection{}
It is well known that every finite spherical group is isomorphic (as an abstract group) to a member of the four infinite families 
\qtnl{221022b}
C_n,\quad D_{2n},\quad C_n\times C_2,\quad D_{2n}\times C_2\qquad (n\ge 1),
\eqtn
or to one of the six groups
\qtnl{221022c}
\alt(4),\quad \sym(4),\quad \alt(5),\quad \alt(4)\times C_2,\quad \sym(4)\times C_2,\quad \alt(5)\times C_2.
\eqtn
A complete classification of finite spherical groups as permutation groups was obtained  in~\cite{Klavik2022}. The lemma below is just a reformulation of this result in a form convenient for subsequent application.

\lmml{221022a}
Let $G$ be a finite spherical noncyclic group and $\Delta\in\orb(G)$. Assume that $\Delta$ is not a  faithful regular orbit. 
\nmrt
\tm{1}  If $G\cong  C_n$ or $G\cong  C_n\times C_2$, then the group $G^\Delta$ is cyclic.
\tm{2}  If $G\cong  D_{2n}$, $n\ge 3$, then $|\Delta|=1$, $2$ or $n$, and $\Delta$ is faithful if $|\Delta|=n$.
\tm{3}  If  $G$ is not as in (1) and (2), then the number $|\Delta|$ and the stabilizer $G_\delta$, $\delta\in\Delta$, are as in the fourth column of Table~\ref{fig:table}.
\enmrt
\elmm

The proof of Theorem~\ref{310722b} is based on the three theorems below and is presented in the end of the section. The first of them states that the coherent configuration~$\cX$ associated with a spherical group $G$ belonging to one of the families~\eqref{221022b} or isomorphic to $\alt(4)\times C_2$ is separable (Theorem~\ref{221022x}). Then we prove that if $G$ is isomorphic (as abstract group) to some groups in~\eqref{221022c}, then $G$ is strongly spherical only if it has a faithful regular orbit (Theorem~\ref{221022y}).  Finally, we prove that the domination free groups among all  remaining groups $G$ are strongly spherical and the coherent configuration $\cX$ has an $\mS_2$-representation which is rigid and faithful (Theorem~\ref{221022z}). 

\begin{table}
\begin{center}
{\tiny
\begin{tabular}[t]{|c|c|l|c|}
\hline
$G$      &   $|G|$      & generators and relations of $G$ &  orbit lengths / point stabilizers \\
\hline
$D_{2n}\times C_2$    & $4n$    &  $\langle t_1,t_2,t_3:\ t_1^2=t_2^2=t_3^2=$   &  $2n^\infty$, $2n$, $2n$, $n$, $n$, $2$\\ 
&   & $(t_1t_2)^n=(t_1t_3)^2=(t_2t_3)^2=1\rangle$     &  $\grp{t_1}$,  $\grp{t_2}$,  $\grp{t_3}$,   $\grp{t_2,t_3}$,  $\grp{t_1,t_3}$, $\grp{t_1,t_2}$\\
\hline
$\alt(4)$    &  $12$     &  $\langle r_1,r_2:\ r_1^3=r_2^3=$  & $6$, $4^2$\\ 
&   &  $(r_1r_2)^2=1\rangle$   &  $\grp{r_1r_2}$,  $\grp{r_1}$\\
\hline
$\sym(4)$    &  $24$    & I:  $\langle t_1,t_2,t_3:\ t_1^2=t_2^2=t_3^2=$   & $12^\infty$, $6$, $4^2$\\ 
&   & $(t_1t_3)^3=(t_1t_2)^3=(t_2t_3)^2=1\rangle$    &  $\grp{t_1}$,  $\grp{t_2,t_3}$,  $\grp{t_1,t_2}$ \\
\hhline{|~| ~|-|-|}
&      & II:  $\langle r_1,r_2:\ r_1^3=r_2^4=$   & $12$, $8$, $6$\\ 
&   & $(r_1r_2)^2=1\rangle$  &  $\grp{r_1r_2}$,  $\grp{r_1}$,  $\grp{r_2}$\\
\hline
$\alt(5)$    &  $60$     &  $\langle r_1,r_2:\ r_1^3=r_2^5=$  & $30$, $20$, $12$\\ 
&   &  $(r_1r_2)^2=1\rangle$   &  $\grp{r_1r_2}$,  $\grp{r_1}$, $\grp{r_2}$\\
\hline
$\alt(4)\times C_2$    &  $24$     &  $\langle r,t,z:\ r^3=t^2=z^2=$  & $12^\infty$, $8$, $6$\\ 
&   &  $(rt)^3=[z,r]=[z,t]=1\rangle$   &  $\grp{z}$,  $\grp{r}$, $\grp{z,t}$\\
\hline
$\sym(4)\times C_2$    &  $48$     &  $\langle t_1,t_2,t_3:\ t_1^2=t_2^2=t_3^2=$  & $24^\infty$, $24^\infty$, $12$, $8$, $6$\\ 
&   & $(t_1t_2)^4=(t_1t_3)^3=(t_2t_3)^2=1\rangle$     &  $\grp{t_1}$,  $\grp{t_2}$, $\grp{t_2,t_3}$,  $\grp{t_1,t_3}$,   $\grp{t_1,t_2}$\\
\hline
$\alt(5)\times C_2 $    &  $120$     &  $\langle t_1,t_2,t_3:\ t_1^2=t_2^2=t_3^2=$  & $60^\infty$, $30$, $20$, $12$\\ 
&   & $(t_1t_2)^5=(t_2t_3)^2=(t_1t_3)^3=1\rangle$     &  $\grp{t_1}$,  $\grp{t_2,t_3}$,  $\grp{t_1,t_3}$,   $\grp{t_1,t_2}$\\
\hline
\end{tabular}
}	
\end{center}
\caption{Orbits of a spherical permutation group $G$ (the superscript on a number in the fourth column indicates the maximal number of $G$-orbits equivariant to a given one; the index $1$ is omitted).}
\label{fig:table}
\end{table}

\subsection{} In what follows, $G$ is  a spherical permutation group on~$\Omega$, and $\cX=\inv(G)$ is the coherent configuration associated with~$G$.

\thrml{221022x}
Let $G$ belong to one of the families~\eqref{221022b} or $G\cong\alt(4)\times C_2$. Then the coherent configuration $\cX$ is separable.
\ethrm
\prf Let $G\cong C_n$ or $G\cong C_n\times C_2$.\footnote{This includes the case $G=1$.} By Lemma~\ref{181022a}, we may assume that $G$ has no faithful regular orbits. By Lemma~\ref{221022a}(1), this implies that every transitive constituent of $G$ is a regular cyclic group. Then  the lattice of the normal subgroups of the constituent is distributive. Hence the coherent configuration $\cX$ satisfies the hypothesis of~\cite[Theorem~1.2]{Hirasaka2018}, whence it follows that $\cX$ is separable. 

Let $G\cong D_{2n}$, $n\ge 3$.   By Corollary~\ref{200922a}, we may assume that $G$ has no orbits of cardinalities~$1$ and~$2$. By Lemma~\ref{221022a}(2),  this implies that every orbit of~$G$ is faithful of cardinality~$n$.  Therefore,  the coherent configuration $\cX$ is separable by Lemma~\ref{240922b}.

Let $G\cong D_{2n}\times C_2$, $n\ge 2$. By Lemma~\ref{310722d} and Corollary~\ref{200922a}, we may assume that~$G$ is domination free and has no orbits of cardinality~$2$. Then by Lemma~\ref{221022a}(3), we have $|\orb(G)|\le 5$ and  
$$
\orb(G)\subseteq\{\Delta_1,\Delta_2,\Delta_3,\Gamma_1,\Gamma_2\},
$$
where the notation is chosen so that  the point stabilizers of the transitive constituents $G^{\Delta_i}$, $1\le i\le 3$, and $G^{\Gamma_j}$, $1\le j\le 2$,  are equal to $\grp{t_i}$ and $\grp{t_j,t_3}$, respectively (see the first row of Table~\ref{fig:table}). Note that $\Delta_1=\Delta_2$ and $\Gamma_1=\Gamma_2$ if $n$ is odd. 

First, assume that neither $\Delta_1,\Delta_2\not\in\orb(G)$. Since $G$ is domination free and  $\Gamma_1$ and $\Gamma_2$ are dominated by $\Delta_3$, it follows that if $\Delta_3\in\orb(G)$, then  $\Gamma_1,\Gamma_2\not\in\orb(G)$. Hence $\Omega=\Delta_3$ and the group $G$ is regular. Thus the coherent configuration $\cX=\inv(G,\Delta_3)$ is separable by Lemma~\ref{181022a}. On the other hand, if $\Delta_3\not\in\orb(G)$, then  $\Omega=\Gamma_1\cup\Gamma_2$ and~$G$ is  a permutation group isomorphic to $D_{2n}$.  Thus, the coherent configuration $\cX=\inv(G,\Gamma_1\cup\Gamma_2)$  is separable by Lemma~\ref{240922b}.

Now let $\Delta:=\Omega\cap(\Delta_1\cup\Delta_2)$ be nonempty. Let $\sigma\in G$ be the permutation of~$\Omega$, corresponding to the generator~$t_3$. Clearly, $\sigma$ is a central involution of~$\aut(\cX)$, which acts trivially outside~$\Delta$. Put $e=e(\sigma)$, see the notation before Lemma~\ref{231022a}.\medskip

{\bf Claim 1.} {\it The coherent configuration $\cX_{\Omega/e}$ is separable.}
\prf
Let $\pi$ be the surjection~\eqref{231022w}. Then $\pi$ is identical outside  $\Delta$, and $\pi(\Delta_i)$ and $\pi(\Gamma_j)$ are orbits of the permutation group $\bar G\cong G/\grp{\sigma}$ induced by the natural action of~$G$ on~$\Omega/e$. Inspecting the point stabilizers of~$\bar G$, it is not difficult to verify that  for $i=1,2$,
\nmrt
\tm{i} if $\Delta_i,\Gamma_i\subseteq\Omega$, then  $\pi(\Delta_i)$ and $\pi(\Gamma_i)$ are
equivariant,
\tm{2} if $\Delta_3,\Gamma_i\subseteq\Omega$, then $\pi(\Gamma_i)$ is dominated by $\pi(\Delta_3)$.
\enmrt
Now if $\Delta_3\subseteq\Omega$, then $\pi(\Delta_3)$ is a faithful regular orbit of the group $\bar G$. Hence the coherent configuration  $\cX_{\Omega/e}=\inv(\bar G)$  is separable by Lemma~\ref{181022a}. On the other hand, if $\Delta_3\not\subseteq\Omega$, then $\bar G\cong D_{2n}$ and  the coherent configuration  $\cX_{\Omega/e}$  is separable by Lemma~\ref{240922b}.
\eprf

{\bf Claim 2.} {\it The parabolic $e$ has a complement.}
\prf
We assume that $\Delta=\Delta_1\cup\Delta_2$; the cases $\Delta=\Delta_1$ and $\Delta=\Delta_2$ are considered in a similar way. Let $\Lambda$ be the disjoint union of two copies $\Lambda_1$ and~$\Lambda_2$ of the group~$G$, and let $H\le\sym(\Lambda)$ be the (semiregular) group induced by the action of~$G$ on $\Lambda$ by right multiplications. Thus, $\orb(H)=\{\Lambda_1,\Lambda_2\}$ and the transitive constituents of~$H$ are regular. 

Denote by $e_0$ the equivalence relation on~$\Lambda$, the classes of which are the right cosets of~$\grp{t_1}$ in~$\Lambda_1$ and the right cosets of~$\grp{t_2}$ in~$\Lambda_2$. Then $e_0$ is a parabolic of the coherent configuration $\inv(H)$, and identifying in a natural way the sets $\Delta$ and $\Lambda/e_0$, we have
$$
\cX_\Delta=\inv(H)_{\Lambda/e_0}.
$$
The cosets of $\grp{t_3}$ in~$\Lambda_1$ and in~$\Lambda_2$ are the classes of an equivalence relation $e_1$ on~$\Lambda$  that is again a parabolic of~$\inv(H)$. Moreover, if $\pi_0:\Lambda\to\Lambda/e_0$, $\lambda\mapsto \lambda e_0$ is the natural surjection, then $\pi_0(e_1)=e$.

Denote by $\Lambda'_1$ the union of the group $\grp{t_1,t_2}\le G$ contained in $\Lambda_1$ and the same group contained in~$\Lambda_2$. Denote by $\Lambda_2'$ the set theoretical complement  of $\Lambda'_1$ in~$\Lambda$.  This complement is equal to the union of the coset $\grp{t_1,t_2}t_3$ contained in $\Lambda_1$ and the same coset in~$\Lambda_2$.  In these notations, $\Lambda$ is split into the disjoint union of the two sets $\Lambda'_1$ and $\Lambda'_2$, and the equivalence relation $e_2$ corresponding to this partition is a parabolic of~$\inv(H)$. It is easily seen that
$$
e_1\cap e_2=e_0
$$
implying that the intersection of any class of $e_2$ with any class of $e_1$ is a class of $e_0$. Thus, the parabolic $\pi_0(e_2)$ is a complement of the parabolic $\pi_0(e_1)=e$.
\eprf

By Claims 1 and 2, the coherent configuration $\cX$ satisfies the condition of Lemma~\ref{231022a}, and hence is separable, which completes the case $G\cong D_{2n}\times C_2$.

Finally, let $G\cong\alt(4)\times C_2$. By Lemma~\ref{310722d} and Corollary~\ref{200922a}, we may assume that $G$ is domination free. By Lemma~\ref{221022a}(3), this implies that  $|\orb(G)|\le 3$ and the set 
$$
N(G)=\{|\Delta|:\ \Delta\in\orb(G)\}
$$
 is contained in the set~$\{6,8,12\}$. The rest of the proof is based on computations in the {\sf GAP}-package COCO2P~\cite{KlinCOCO2P}. First, we have
$$
\aut(\cX)\cong\css
\alt(4)                        &\text{if } N(G)\subseteq\{6,12\},\\
\sym(4)\times C_2  &\text{if } N(G)=\{8\}.\\
\ecss
$$
This leaves us with three  cases: $N(G)=\{6,8\}$, $\{8,12\}$, and $\{6,8,12\}$.  Furthermore, the orbit of cardinality $6$ is dominated by the orbit of cardinality~$12$. Thus either $N(G)=\{6,8\}$ or $\{8,12\}$. 

In each of the two cases, the group $\aut(\cX)$ has a unique central involution~$\sigma$, and  (in the notation of Lemma~\ref{231022a}) $|\Delta|=8$ and the parabolic~$e=e(\sigma)$ of $\cX$  has a complement.  Finally,  the quotient coherent configuration  $\cX_{\Omega/e}$ is separable: if  $N(G)=\{8,12\}$, then this follows from Lemma~\ref{181022a} (the group  $\aut(\cX_{\Omega/e})$ has faithful regular orbit), whereas if $N(G)=\{6,8\}$, then this follows by inspecting the database of small coherent configurations~\cite{Ziv-av2018}. Thus the coherent configuration~$\cX$ satisfies the condition of Lemma~\ref{231022a}, and hence is separable. This completes the  proof of Theorem~\ref{221022x}.
\eprf

\subsection{} 
The proof of  Theorems~\ref{221022y} and~\ref{221022z} is based on case-by-case computer computations. They use  the  {\sf GAP} packages COCO2P~\cite{KlinCOCO2P} for computations  with coherent configurations and WLFast which is a faster implementation of the Weisfeiler-Lehman algorithm provided to us by S.~V.~Skresanov, and  computer systems Mathematica~\cite{Mathematica} and Maple~\cite{Maple2022}. The protocols of the computations can be found at~\url{https://zemanpeter.github.io/assets/files/comp.zip}.

\thrml{221022y}
Let $G\cong\alt(4)$, $\sym(4)$ in the representation~II, or $\alt(5)$. Assume that $G$ has no faithful regular orbit. Then $G$ is not strongly spherical unless $G$ is induced by the action of~$\sym(4)$ of degree~$12$.\footnote{The permutation groups induced   by the action of~$\sym(4)$ of degree~$12$  in the representations~I and~II are permutation isomorphic.}
\ethrm
\prf
The assumption of the theorem implies by Lemma~\ref{221022a}(3) that $G$ has at most three orbits. All possible cases are given in Table~\ref{tbl10}, the second column of which represents the set $N(G)$. The computation starts with constructing  the coherent configuration $\cX=\inv(G)$; the rank and the automorphism group of $\cX$ are indicated in the third and fourth columns, respectively. Note that $G$ is strongly spherical only if 
\qtnl{261022a}
G=\aut(\cX).
\eqtn 
Therefore comparing $G$ and $\aut(\cX)$ reduces the number of possibilities for~$G$, e.g., the degree of a strongly spherical permutation group isomorphic to $\alt(4)$ can be equal to $10$ or $20$ only (see the first five rows of Table~\ref{tbl10}).
\begin{table}[t]
\begin{center}
{\footnotesize
\begin{tabular}[t]{|c|c|c|c|c|c|}
\hline
type             &  $\orb(G)$      &  $\rk(\cX)$    & $\aut(\cX)$                   & $|\fX'(G)|$ &  $|\fX(G)|$    \\
\hline
$\alt(4)$    & $6$                  & $4$                 &   $C_2\times \alt(4)$     &                     & \\
\hhline{|~|-|-|-|--|}
& $4$                  & $2$                  &   $\sym(4)$                    &                     &   \\
\hhline{|~|-|-|-|--|}
& $4+4$             & $8$                  &   $\sym(4)$                     &                     &   \\
\hhline{|~|-|-|-|--|}
& $6+4$              & $10$                &   $\alt(4)$                       &  $19$          & $0$    \\
\hhline{|~|-|-|-|--|}
& $6+4+4$         & $20$                 &   $\alt(4)$                      &  $294$       &$0$   \\
\hline
$\sym(4)$II   & $6$              & $3$                   &   $C_2\times \sym(4)$ &       & \\
\hhline{|~|-|-|-|--|}
& $8$                   & $4$                &   $C_2\times \sym(4)$ &       &  \\
\hhline{|~|-|-|-|--|}
& $12+8$            & $19$               &   $ \sym(4)$                   &   $330$  & $14$ \\
\hhline{|~|-|-|-|--|}
& $12+6$            & $16$                &   $ \sym(4)$                    &  $148$  &  $6$ \\
\hhline{|~|-|-|-|--|}
& $8+6$                   & $11$          &   $C_2\times \sym(4)$ &            &  \\
\hhline{|~|-|-|-|--|}
& $12+8+6$            & $32$          &   $ \sym(4)$                 & $4827$   &  $36$ \\
\hline
$\alt(5)$      & $30$                     &  $16$          &  $\alt(5)$                     & $52$        & $5$  \\
\hhline{|~|-|-|-|--|}
&$20$                      &   $8$           &  $\alt(5)$                     &   $10$       &  $0$  \\
\hhline{|~|-|-|-|--|}
& $12$                     &    $4$          &  $C_2\times\alt(5)$    &                   &  \\
\hhline{|~|-|-|-|--|}
& $30+20$              &   $44$         &  $\alt(5)$                     & $2385$      &  $148$   \\
\hhline{|~|-|-|-|--|}
& $30+12$            &  $32$           &   $\alt(5)$                    & $601$        & $25$  \\
\hhline{|~|-|-|-|--|}
& $20+12$              &  $20$          &    $\alt(5)$                  & $80$            &   $1$ \\
\hhline{|~|-|-|-|--|}
& $30+20+12$       &   $68$          &   $\alt(5)$                  &  $46631$       &  $685$ \\
\hline
\end{tabular}
}	
\end{center}
\caption{The groups $G\cong\alt(4)$, $\sym(4)$II, or $\alt(5)$.}
\label{tbl10}
\end{table}

For each remaining group $G$, we need to verify that there is no polyhedral graph~$X$ with vertex set~$\Omega$, such that $\WL(X)=\cX$.  Note that every such $X$ satisfies the following conditions:
$$
E\in S^\cup, \quad\WL(X)=\cX,\quad |E|\le  3|\Omega|-6,
$$
where $S$ is the set of basis relations of~$\cX$. The first condition allows us to search for possible $E$ among the subsets of~$S$. Since every subgraph of $X$ is planar, we can strengthen the third  inequality by requiring that
$$
|E\cap(\Delta\times\Gamma)|\le |3|\Delta|-6
$$
for all $\Delta,\Gamma\in\orb(G)$. Denote $\fX'=\fX'(G)$ the set of all such graphs~$X$. 

The action of the group $\Phi$ of algebraic automorphisms of the coherent configuration~$\cX$ on $S^\cup$ induces a natural action of this group on the set~$\fX'$. Let $\fX(G)$ be the set of distinct representatives~$X$ of $\orb(\Phi,\fX')$, such that $\WL(X)=\cX$. Since every algebraic isomorphism of $\cX$ is induced by isomorphism (this can also be checked with the help of COCO2P), we conclude that the group~$G$  is strongly spherical if and only if the set $\fX(G)$ contains a polyhedral graph.

To complete the proof, we compute the set $\fX(G)$ for all groups $G$ in Table~\ref{tbl10}, that satisfy~\eqref{261022a}; the number $|\fX(G)|$ is given in the fifth column of Table~\ref{tbl10}. Then using the functions ``IsPlanar'' and ``VertexConnectivity'' included in  Maple, we verified that  $\fX(G)$ contains no polyhedral graphs.
\eprf

\subsection{} In this subsection, we consider the remaining groups~$G$ and prove Theorem~\ref{221022z}. A key role here is played  by the concept of $\mS_2$-representation studied in Section~\ref{290222a}.

\thrml{221022z}
Let $G\cong\sym(4)$ in the representation~I, $\sym(4)\times C_2$, or $\alt(5)\times C_2$. Assume that $G$ is domination free. Then $G$ is strongly spherical and the coherent configuration $\cX$ has an $\mS_2$-representation which is rigid and faithful.
\ethrm
\prf
The assumption of the theorem implies by Lemma~\ref{221022a}(3) that $G$ has at most three orbits. All possible cases are given in Table~\ref{tbl4}, the second column of which represents the set $N(G)$.  
\begin{table}[t]
\begin{center}
{\tiny%\scriptsize %\footnotesize
\begin{tabular}{|c|c|c|l|c|c|c|}
\hline
type & $\orb(G)$      &  $|S| / |S_\rho|$      &  $X$  &  $\deg(X_\rho)$ &  $A_\rho$ & $\lambda_\rho$    \\
\hline
$\sym(4)$I    & $12$     &   $7/6$     &  Truncated Tetrahedral  &  $12[3]$           & L   & $1$  \\
\hhline{|~|-|-|-|-|-|-|}
& $4$       &   $2/2$     &  Tetrahedral                     &  $4[3]$            & L  & $4$ \\
\hhline{|~|-|-|-|-|-|-|}
& $6+4$    &   $9/9$     & $X_{10}$, polyhedral    &  $6[6]+4[3]$  & A  & $\sqrt{2}$ \\
\hline
$\sym(4)\times$  & $24$I   &   $14/12$  &  Small Rhombicuboctah.  & $24[4]$  & L & $1$ \\
\hhline{|~|-|-|-|-|-|-|}
$C_2$                   & $24$II   &   $16/11$  &  Truncated Octahedral    & $24[3]$          & L   & $4-\sqrt{2}$  \\
\hhline{|~|-|-|-|-|-|-|}
& $12$       &  $5/5$      & Cuboctahedral                 & $12[4]$           & L   & $2$  \\
\hhline{|~|-|-|-|-|-|-|}
& $8$         &  $4/4$      & Cubical                            &  $8[3]$            & L   &  $2$ \\
\hhline{|~|-|-|-|-|-|-|}
& $6$         &  $3/3$       & Octahedral                     &   $6[4]$           &  L   & $4$ \\           
\hhline{|~|-|-|-|-|-|-|}
& $24$I$+24$II  & $54/49$  & $X_{48},$ polyhedral & $24[6]+24[4]$ & $A$ & $(3+\sqrt{33})/2$ \\
\hhline{|~|-|-|-|-|-|-|}
& $24$II$+8$     & $28/25$  &  $X_{32},$ polyhedral & $24[5]+8[6]$  &  A & $1-\sqrt{10}$\\
\hhline{|~|-|-|-|-|-|-|}
& $12+8$             & $15/14$   & $X_{20},$ polyhedral & $12[6]+ 8[3]$  & A  & $1+\sqrt{5}$ \\
\hhline{|~|-|-|-|-|-|-|}
& $12+6$              &  $14/12$  & $X_{18},$ polyhedral & $12[6]+ 6[4]$   & A & $1-\sqrt{5}$ \\
\hhline{|~|-|-|-|-|-|-|}
& $8+6$                &  $11/11$   & Rhombic Dodecahedral & $8[3]+ 6[4]$ & L &  $(7-\sqrt{17})/2$ \\
\hhline{|~|-|-|-|-|-|-|}
& $12+8+6$         &  $28/25$  & Disdyakis Dodecahedral & $12[4]+8[6]+6[8]$  & A & $4$ \\
\hline
$\alt(5)\times$  & $60$  &   $32/22$  &  Small Rhombicosidodecah.   &  $60[4]$  & L & $(3-\sqrt{5})/2$ \\
\hhline{|~|-|-|-|-|-|-|}
$C_2$                 & $30$  &   $10/9$   &  Icosidodecahedral                  &  $30[4]$  & L  &  $3-\sqrt{5}$ \\
\hhline{|~|-|-|-|-|-|-|}
& $20$  &   $6/6$     &   Dodecahedral                          &  $20[3]$  &L  & $3-\sqrt{5}$ \\
\hhline{|~|-|-|-|-|-|-|}
& $12$  &    $4/4$     &  Icosahedral                              &  $12[5]$  &L  & $5-\sqrt{5}$ \\
\hhline{|~|-|-|-|-|-|-|}
& $30+20$  &   $30/28$  &  $X_{50},$ nonplanar   &  $30[2]+20[3]$  & A & $(\sqrt{10}-\sqrt{2})/2$  \\
\hhline{|~|-|-|-|-|-|-|}
& $30+12$  &   $24/22$  &   $X_{42},$ planar    &  $30[2]+12[5]$  & A & $\sqrt{5-\sqrt{5}}$  \\
\hhline{|~|-|-|-|-|-|-|}
& $20+12$  &   $18/18$  &  Rhombic Triacotahedral  &  $20[3]+12[5]$  & L & $3-\sqrt{5}$ \\
\hhline{|~|-|-|-|-|-|-|}
& $30+20+12$  &  $52/49$  &   Disdyakis Triacontahedral   &$30[4]+20[6]+12[10]$  &A &   $\sqrt{5}$ \\
\hline
\end{tabular}
}
\end{center}
\caption{The groups  $G\cong\sym(4)$I, $\sym(4)\times C_2$, $\alt(5)\times C_2$ }
\label{tbl4}
\end{table}

\begin{figure}[t!]
\centering
\includegraphics[scale=0.35]{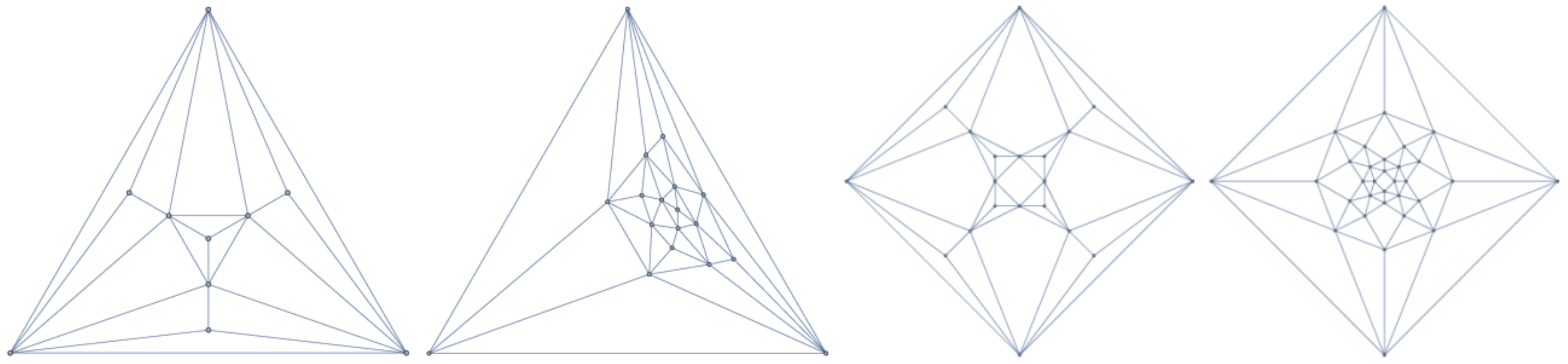}
\\
\includegraphics[scale=0.40]{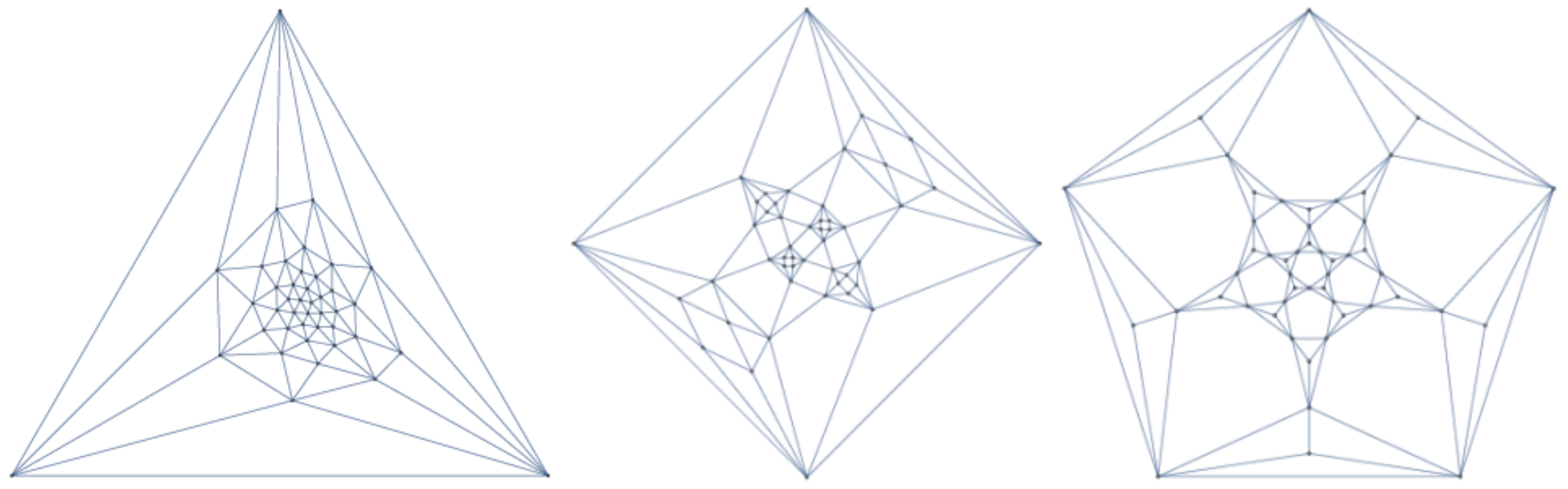}
\caption{The polyhedral  graphs  $X_{10},X_{18}, X_{20},X_{32},Y_{42},X_{48}, Y_{50}$.}
\label{X4250}
\end{figure}

To prove that the group $G$ is strongly spherical, one needs to find a polyhedral graph $X$ such that the coherent configurations $\WL(X)$ and $\cX$ are isomorphic. Except for the seven cases, $X$ can be found among the known polyhedral graphs  the commonly used names of which are given in the fourth column of Table~\ref{tbl4}; for the other seven cases we take the graphs shown in Fig.~\ref{X4250} (the indices at $X$ and~$Y$ are equal to the number of vertices of the graph). 

Using the Mathematica package, we first verified that each  graph in fourth column of Table~\ref{tbl4} is polyhedral\footnote{Certainly, this is not necessary to do for the known graphs.} and construct its the adjacency matrix $A=A(X)$. Then using the COCO2P function ColorGraphByMatrix and the WLFast function WLStabilizationByHash, we construct the coherent configuration $\WL(X)$. Finally, we verified that  $\WL(X)\cong\cX$ with the help of the COCO2P function IsColorIsomorphicColorGraph. This completes the proof of the first statement.

\begin{figure}[h]
\centering
\includegraphics[scale=0.35]{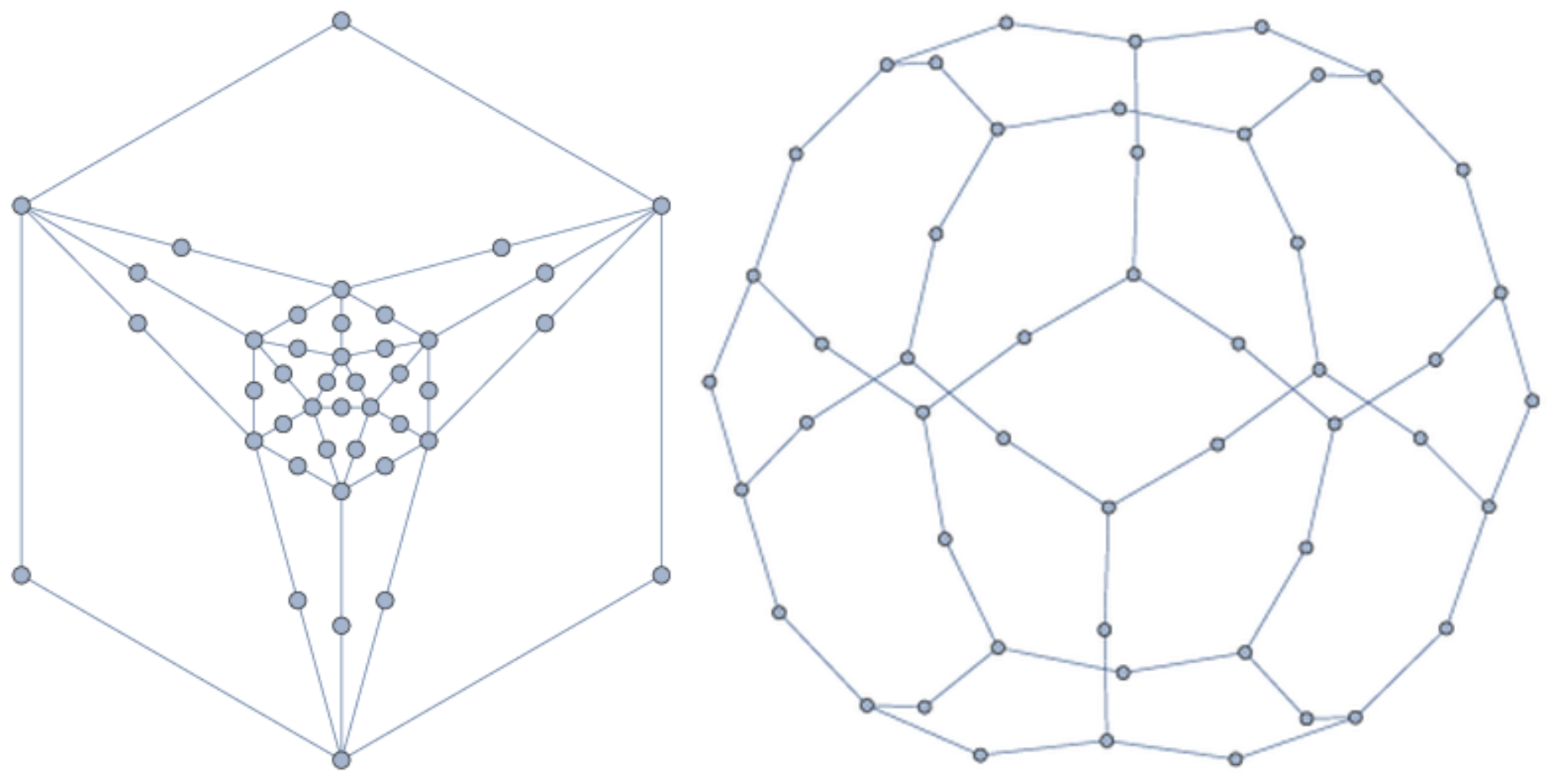}	
\caption{The nonpolyhedral  graphs  $X_{42}, X_{50}$.}
\label{XY}
\end{figure}

To prove the second statement, we make use of Lemma~\ref{030922a} to construct $\mS_2$-representation of the coherent configuration~$\cX$. As the matrix $P$, we take the normalized orthogonal projection $P_\rho$ to the eigenspace of a matrix $A_\rho$ which is the Laplacian (``L'' in the sixth column of Table~\ref{tbl4}) or adjacency matrix (``A'' in the sixth column of Table~\ref{tbl4}) of the graph~$X$, corresponding to an eigenvalue $\lambda_\rho$ (the seventh columns of Table~\ref{tbl4}) of multiplicity~$3$. For the polyhedral graphs~$Y_{42}$ and~$Y_{50}$, the computation of $P_\rho$ was not successful.\footnote{The symbolic linear algebra in the Maple package  was unable to calculate the projection due to too much calculation.} Because of this, these graphs were replaced by the graphs~$X_{42}$ and~$X_{50}$, see Fig.~\ref{XY}. For each of them, $\WL(X)$ and $\cX$ are still isomorphic,  and the adjacency matrix still has an eigenvalue $\lambda_\rho$ of multiplicity~$3$.
The distributions of vertex degrees of all the graphs~$X$ are indicated in the fifth column of Table~\ref{tbl4}.

The remaining calculations are organized as follows. In view of the isomorphism $\WL(X)\cong\cX$, we may assume that $\cX=\WL(X)$. Inspecting the entries of the matrix $P_\rho$ (they are algebraic numbers and therefore we need symbolic linear algebra in Maple), we compute the relations~\eqref{301022s} in the form of unions of some basis relations of~$\cX$, and then  the rainbow $\cX_\rho\le\cX$, see Section~\ref{290222a}. The ranks of $\cX$ and $\cX_\rho$ are presented in the third column of Table~\ref{tbl4}. If these ranks coincide, then $\cX=\cX_\rho$ and the $\mS_2$-representation~$\rho$ of the coherent configuration~$\cX$, corresponding to~$P_\rho$, is obviously faithful; otherwise, we compute  $\WL(\cX_\rho)$ by using the package  WLFast and test that in all cases,  $\WL(\cX_\rho)=\WL(X)$. Thus the $\mS_2$-representation $\rho$ is faithful.

Finally, to verify that $\rho$ is rigid, denote by $T_\rho$ the set of all $s\in S$ contained in the antipodal parabolic with respect to~$\rho$, and by $\Pi_\rho$ the partition of~$S$ such that $r$ and $s$ belong to the same class if and only if $w(r)=w(s)$. We defined a new COCO2P function which given a coherent configuration $\cX$, a relation $s\in S$, a partition $\Pi_\rho$, and the set $T_\rho$, tests whether a set $\Delta=\{\alpha,\beta\}$ for some $(\alpha,\beta)\in s$ is $\rho$-rigid (see Section~\ref{290222a}). The computations shows that for each graph~$X$, there is at least one $s\in S$, for which $\Delta$ is $\rho$-rigid. This shows that  the $\mS_2$-representation $\rho$ is rigid. This completes the proof.
\eprf

\subsection{Proof of Theorem~\ref{310722b}.}
Let $G$ be a strongly spherical permutation group and $\cX=\inv(G)$. Then  $G$ belongs to one of the families~\eqref{221022b}, or is isomorphic to one of the groups~\eqref{221022c}. By Theorem~\ref{221022x}, we may assume that neither $G$ belong to one of the families~\eqref{221022b} nor $G\cong\alt(4)\times C_2$. 

Now if  $G\cong\alt(4)$, $\sym(4)$ in the representation~II, or $\alt(5)$, then $G$ has a faithful regular orbit by Theorem~\ref{221022y} and hence the coherent configuration~$\cX$ is separable by Lemma~\ref{181022a}. 

Finally, in the remaining case, we may assume by Lemma~\ref{310722d} that $G$ is domination free. By Theorem~\ref{221022z}, this implies that the coherent configuration $\cX$ has an $\mS_2$-representation which is both rigid and faithful. Thus, $\cX$ is separable by Theorem~\ref{070922v}. {\hfill$\square$}

%\bibliography{bibinp}{}
%\bibliographystyle{amsplain}

\end{document}